\documentclass[12pt,notitlepage]{amsart}

\usepackage{latexsym,amsfonts,amssymb,amsmath,amsthm} 
\usepackage{graphicx}
\usepackage{amsmath}
\usepackage{amssymb}
\usepackage{mathrsfs}

\newtheorem{thm}{Theorem}[section]

\newtheorem{prop}[thm]{Proposition}
\newtheorem{lemma}[thm]{Lemma}

\newtheorem*{thm*}{Theorem}

\newcommand{\lam}{\lambda}

\begin{document}
\title{Bounding $|\zeta(\tfrac 12+it)|$ on the Riemann hypothesis}

\author{Vorrapan Chandee} 
\address{Department of Mathematics \\ Stanford University\\
450 Serra Mall, Bldg. 380\\
Stanford, CA 94305-2125}
 \email{vchandee@math.stanford.edu}
\author{K. Soundararajan}
\address{Department of Mathematics \\ Stanford University\\
450 Serra Mall, Bldg. 380\\
Stanford, CA 94305-2125}
 \email{ksound@math.stanford.edu}
\thanks
{The authors are partially supported by a grant from the NSF  (DMS-0500711)}
\maketitle
\section{Introduction}

\noindent In 1924  Littlewood \cite{Littlewood} proved that the Riemann Hypothesis (RH) implies a strong 
form of the Lindel{\" o}f hypothesis; namely, on RH, for large $t$ there is a constant $C$ 
such that 
\begin{equation} 
\label{LitLH}
|\zeta(\tfrac 12+it)| \ll \exp\Big( C\frac{\log t}{\log \log t}\Big). 
\end{equation} 
In the intervening years no improvement has been made over \eqref{LitLH}, 
except in reducing the permissible value of $C$, see \cite{RamaSan, Moment}.   
In \cite{Moment} Soundararajan showed that \eqref{LitLH} holds for any $C> (1+\lambda_0)/4 = 0.372\ldots$ 
where $\lam_0 =0.4912\ldots$ is the unique positive real number satisfying $e^{-\lam_0} = \lam_0 +\lam_0^2/2$.    In \cite{Fai} Chandee has provided an explicit version of this bound 
for general $L$-functions.
 
A similar situation exists for $S(t) = \frac 1{\pi} \text{arg} \zeta(\frac 12+it)$, where 
the argument is defined by continuous variation along the line segments joining $2$, $2+it$, 
and $\frac 12+it$, taking the argument of $\zeta(s)$ at $2$ to be zero.  On RH Littlewood showed 
that   $S(t) \ll \log t/\log \log t$, and again this bound has not 
been improved except for the size of the implied constant.  Recently Goldston and Gonek \cite{GG}
gave an elegant argument leading to the bound $|S(t)| \le (\frac 12+o(1)) \log t/\log \log t$.  
Their method used the explicit formula together with certain optimal majorants and minorants of characteristic functions of intervals that were constructed by Selberg.    The Goldston-Gonek
result may reasonably be thought of as having attained the limit of existing methods of bounding $S(t)$, 
although it seems likely that the true maximal size of $S(t)$ is even smaller, perhaps $\ll \sqrt{\log t\log \log t}$ (see \cite{FGH}).  

In \cite{Moment} Soundararajan asked for a corresponding treatment for $|\zeta(\tfrac 12+it)|$ 
which would represent the limit of existing methods for bounding $|\zeta(\tfrac 12+it)|$ on RH.  
In this note we present such an approach.  Using Hadamard's factorization formula 
and the explicit formula, we show how the problem of bounding $|\zeta(\tfrac 12+it)|$ 
may be framed in terms of minorizing the function $\log \frac{4+x^2}{x^2}$  by functions whose Fourier 
transforms are supported in a given interval, and drawing upon recent work of Carneiro 
and Vaaler \cite{CV} we find the optimal such minorant.  

\begin{thm}
\label{zetabound}  Assume RH.
For large real numbers $t$ we have
 $$ 
 |\zeta(\tfrac 12+it)| \ll \exp\Big( \frac{\log 2}{2} \frac{\log t}{\log \log t} + O\Big(\frac{\log t \log \log \log t}{(\log \log t)^2}\Big)\Big). 
 $$ 
 \end{thm}

As with $S(t)$, the true maximal size 
of $|\zeta(\tfrac 12+it)|$ may be much smaller, perhaps of size $\exp( \sqrt{(\tfrac 12+o(1))\log t \log \log t})$ 
as suggested by Farmer, Gonek, and Hughes \cite{FGH}.   On the other hand, it is known that 
there are arbitrarily large $t$ such that $|\zeta(\tfrac 12+it)|  \ge \exp((1+o(1)) \sqrt{\log t/\log \log t})$, see \cite{S2}.

 \section{Proof of Theorem \ref{zetabound}}
 
 \noindent Let $\xi(s) = \tfrac 12 s(s-1) \pi^{-s/2} \Gamma(\tfrac s2) \zeta(s)$ 
 denote Riemann's $\xi$-function which is entire of order $1$, satisfies the 
 functional equation $\xi(s)=\xi(1-s)$, and whose zeros 
 are the non-trivial zeros of $\zeta(s)$.   Recall (see Chapter 12 of \cite{D}, for example)
 Hadamard's factorization formula 
 \begin{equation*}
 \label{Had1} 
 \xi(s) = e^{A+Bs} \prod_{\rho} \Big(1-\frac{s}{\rho}\Big)e^{s/\rho},
 \end{equation*} 
 where $\rho$ runs over the non-trivial zeros of $\zeta(s)$, and $B= - \sum_{\rho} \text{Re} (1/\rho)$. 
 (Note that Re$(1/\rho)$ is positive and $\sum_{\rho} \text{Re}(1/\rho)$ converges.) 
 We apply this with $s=\tfrac 12+it$ and $s=-\frac 32+it$ and divide.    The absolute convergence 
 of the product allows us to divide term by term, and we find, writing (on RH) $\rho=\tfrac 12+i\gamma$,
 \begin{equation*} 
\Big| \frac{\xi(\frac 12+it)}{\xi(-\frac 32+it)}\Big| = e^{2B} \prod_{\rho} \Big|\frac{i(\gamma-t)}{2+i(\gamma-t)}\Big| e^{\text{Re}(2/\rho)} 
= \prod_{\rho} \Big| \frac{(t-\gamma)^2}{4+(t-\gamma)^2}\Big|^{\frac 12}.
\end{equation*}
Since $\xi(-\tfrac 32+it) = \xi(\tfrac 52-it)$, and $|\zeta(\tfrac 52-it)| \asymp 1$,  
 we deduce using Stirling's formula
 that 
 \begin{equation} 
 \label{bound1} 
 \log |\zeta(\tfrac 12+it)| = \log t + O(1) - \frac 12 \sum_{\gamma} f(t-\gamma), 
 \end{equation}
 where we have set 
 \begin{equation}
 \label{fdef}
 f(x) = \log \frac{4+x^2}{x^2}. 
 \end{equation}
 
 The proof of Theorem \ref{zetabound} now proceeds by replacing 
 $f(t-\gamma)$ by a carefully chosen function that minorizes it, and 
 then invoking the explicit formula.  The properties of the appropriate 
 minorant function are detailed in the following Proposition 
 which we shall demonstrate in the next section. 
 
 \begin{prop}
 \label{mainprop} Let $\Delta$ denote a positive real number.  There is an entire 
 function $g_{\Delta}$ which satisfies the following properties: 
 
 \noindent (i).  For all real $x$ we have 
 $$
-C \frac{1}{1+x^2} \le g_{\Delta}(x) \le f(x), 
$$ 
for some positive constant $C$.  For any complex number 
 $x+iy$ we have 
 $$
 |g_{\Delta}(x+iy)| \ll  \frac{\Delta^2  }{1+ \Delta |x+iy|}e^{2\pi \Delta|y|}.
 $$
 
 \noindent (ii).  The Fourier transform of $g_{\Delta}$, namely 
 $$
 {\hat g}_{\Delta}(\xi) = \int_{-\infty}^{\infty} g_{\Delta}(x) e^{-2\pi i x \xi } dx,  
 $$ 
 is real valued, equals zero for $|\xi| \ge \Delta$, and satisfies $|{\hat g}_{\Delta}(\xi)| \ll 1$.  
 
 \noindent (iii).  The $L^1$ distance between $g_{\Delta}$ and $f$ equals 
 $$
 \int_{-\infty}^{\infty} (f(x)-g_{\Delta}(x)) dx = \frac{1}{\Delta} \Big(2\log 2 -2 \log (1+e^{-4\pi \Delta})\Big) .
 $$
 \end{prop}
 
 Returning to \eqref{bound1}, we have for any positive $\Delta$ 
 \begin{equation} 
 \label{lower}
 \sum_{\gamma} f(t-\gamma)  \ge \sum_{\gamma} g_{\Delta}(t-\gamma). 
 \end{equation} 
 We now invoke the explicit formula connecting zeros and primes, 
 see Lemma 1 of \cite{GG}, or Theorem 5.12 of \cite{IK} . 
 
 \begin{lemma} \label{lem:explicitformula}
Let $h(s)$ be analytic in the strip $|\text{\rm Im }s| \leq 1/2 + \epsilon$ for some $\epsilon > 0$, 
and such that $|h(s)| \ll (1+|s|)^{-(1+\delta)}$ for 
some $\delta >0$ when $|\text{\rm Re }s|\to \infty$.  Let $h(w)$ be a real-valued for real $w$, and set $\hat{h}(x) = \int_{-\infty}^{\infty} h(w) e^{-2\pi i x w} \> dw$. 
Then 
\begin{eqnarray*}
\sum_{\rho} h(\gamma) &=& h\Big(\frac{1}{2i} \Big) +h\Big(-\frac{1}{2i}\Big) 
- \frac{1}{2\pi} {\hat h}(0) \log \pi  + \frac{1}{2 \pi}  \int_{-\infty}^{\infty} h(u) \text{\rm Re }\frac{\Gamma '}{\Gamma}\Big(\frac{1}{4}   + \frac{iu}{2} \Big) \> du \\
&& - \frac{1}{2\pi}\sum_{n=2}^{\infty} \frac{\Lambda(n)}{\sqrt{n}}\Big(\hat{h}\Big( \frac{\log n}{2\pi }\Big) + \hat{h}\Big( \frac{-\log n}{2\pi }\Big)\Big).
\end{eqnarray*}
\end{lemma} 

 We apply Lemma \ref{lem:explicitformula}, taking $h(z)=g_{\Delta}(t -z)$ so that ${\hat h}(x) 
 = {\hat g}_{\Delta}(-x) e^{-2\pi i xt}$.  From (i) of Proposition \ref{mainprop} 
 we find that $h(1/2i) + h(-1/2i) \ll \Delta^2 e^{\pi \Delta}/(1+\Delta t)$, 
 and using (ii) of Proposition \ref{mainprop} that ${\hat h}(0) \ll 1$.  
 Using Stirling's formula, parts (i) and (iii) of Proposition \ref{mainprop}, and that 
 $\int_{-\infty}^{\infty}f(x)dx = 4\pi$ we have
 \begin{eqnarray*} 
 \frac{1}{2\pi } \int_{-\infty}^{\infty} h(u) \text{Re } \frac{\Gamma^{\prime}}{\Gamma}\Big( 
 \frac{1}{4} + \frac{iu}{2} \Big) \> du 
 &=& \frac{1}{2\pi} \int_{-\infty}^{\infty} g_\Delta(u) (\log t + O(\log (2+|u|)) \> du \\
 &=& 2 \log t - \frac{\log t}{\pi \Delta} \log \Big(\frac{2}{1+e^{-4\pi \Delta}}\Big)+ O(1). 
 \end{eqnarray*}
Using these remarks to evaluate the RHS of \eqref{lower}, and inserting that bound in \eqref{bound1} 
we conclude that 
\begin{equation} 
\label{mainbound}
\log |\zeta(\tfrac 12+it)| \le \frac{\log t}{2\pi \Delta} \log \Big( \frac{2}{1+e^{-4\pi \Delta}}\Big) 
+ \frac{1}{2\pi}\text{Re } \sum_{n=2}^{\infty} \frac{\Lambda(n)}{n^{\frac 12+it}} {\hat g}_{\Delta}
\Big(\frac{\log n}{2\pi} \Big) +O\Big( \frac{\Delta^2 e^{\pi \Delta}}{1+\Delta t} + 1 \Big ). 
\end{equation} 

Since 
 \begin{equation*}
 \sum_{n=2}^{\infty} \frac{\Lambda(n)}{\sqrt{n}} \Big|{\hat g}_{\Delta}\Big(\frac{\log n}{2\pi} 
\Big)  \Big| \ll  e^{\pi \Delta},
 \end{equation*}
 taking $\pi \Delta = \log \log t - 3 \log \log \log t$ in \eqref{mainbound} we 
 obtain our Theorem.

 \section{Proof of Proposition \ref{mainprop}: The work of Carneiro and Vaaler} 
 
 \noindent Given a function from ${\Bbb R}$ to ${\Bbb R}$ Carneiro and Vaaler 
 consider the problem of finding optimal majorants and minorants for this 
 function, with the additional property that the majorants and minorants 
are restrictions to the real axis of complex analytic functions of exponential type 
at most $2\pi$.  The majorants and 
 minorants are to be optimal in the sense of minimizing the $L^1$ distance 
 from the given function.   This problem has a long history, going back to 
 work of  Beurling for the signum function which was rediscovered 
 and used by Selberg to study the case of indicator functions of intervals (see \cite{GV, V, CV}).  
 Carneiro and Vaaler solve the optimization problem for a wide class 
 of functions including our function $f(x)$.  
 
Let $\mu$ be a (non-negative) measure defined on the Borel subsets of 
${\Bbb R}_+$ such that 
\begin{equation}
\label{CV1}
0 < \int_0^{\infty} \frac{\lam}{\lam^2 +1} d\mu(\lam) <\infty. 
\end{equation} 
Let 
\begin{equation*}
\label{CV2} 
f_{\mu}(x) = \int_0^{\infty} (e^{-\lam |x|}-e^{-\lam}) d\mu(\lam),  
\end{equation*}
and define 
\begin{equation*} 
\label{CV3} 
G_{\mu}(z) = \lim_{N\to \infty}\Big( \frac{\cos \pi z}{\pi }\Big)^2  \sum_{n=-N}^{N+1} \Big(\frac{ f_{\mu}(n-\frac 12)}{(z-n+\frac 12)^2} + \frac{f_{\mu}^{\prime}(n-\frac 12)}{(z-n+\frac 12)} \Big) .
\end{equation*}
 Theorem 1.1 of Carneiro and Vaaler then demonstrates that $G_{\mu}(z)$ converges uniformly on 
 compact subsets of ${\Bbb C}$, defines an entire function of exponential type 
 at most $2\pi$, and that for real $x$ we have $G_{\mu}(x) \le f_{\mu}(x)$.  Moreover they 
 show that $G_{\mu}$ minimizes the $L^1$ distance from $f_{\mu}$ (in particular $f_{\mu}-G_{\mu}$ 
 is integrable)
 among all minorants of $f_{\mu}$ with exponential 
 type at most $2\pi$.

 Let $\Delta$ be a given positive real number, and consider the measure 
 $$
 d\mu_{\Delta}(\lam) = \frac{2(1-\cos(2\Delta \lam))}{\lam} d\lam.
 $$
 This measure satisfies \eqref{CV1}, and moreover 
 \begin{equation} 
 \label{Meas1}
 \int_0^{\infty} (e^{-\lam|x|}-e^{-\lam}) \frac{2(1-\cos(2\Delta\lam))}{\lam} d\lam 
 = \log \Big(\frac{4\Delta^2 +x^2}{x^2} \Big)- \log (4\Delta^2 +1). 
 \end{equation}
 The identity \eqref{Meas1} may be checked by noting that both sides equal zero for $x=1$, 
 and that the derivatives of both sides agree (a little care is needed at $x=0$ where the result 
 follows by continuity).  Let us denote the RHS of \eqref{Meas1} by $f_{\Delta}(x) = f(x/\Delta)-f(1/\Delta)$. 
 Let 
 \begin{equation*} 
 \label{Meas2} 
 G_{\Delta}(z) = \Big( \frac{\cos \pi z}{\pi }\Big)^2 \sum_{n=-\infty}^{\infty} 
 \Big( \frac{f_{\Delta}(n-\frac 12)}{(z-n+\frac 12)^2} + \frac{f_{\Delta}^{\prime}(n-\frac 12)}{(z-n+\frac 12)}\Big) 
 \end{equation*} 
 denote the corresponding optimal function of Carneiro and Vaaler.  
   
First we record an upper bound for $G_\Delta(z)$.  By an application of 
the Poisson summation formula we see that 
$$ 
\Big( \frac{\cos \pi z}{\pi}\Big)^2 \sum_{n=-\infty}^{\infty} \frac{1}{(z-n+\frac 12)^2} 
= \sum_{n=-\infty}^{\infty} \Big(\frac{\sin (\pi (z-n+\frac 12))}{\pi(z-n+\frac 12)}\Big)^2 
= 1,
$$
so that 
\begin{equation}
\label{Meas3}
G_{\Delta}(z) + f(1/\Delta) = \sum_{n=-\infty}^{\infty} \Big(\frac{\sin (\pi (z-n+\frac 12))}{\pi(z-n+\frac 12)}\Big)^2 \Big( f\Big(\frac{n-\frac 12}{\Delta} \Big) + \frac{(z-n+\tfrac 12)}{\Delta} f^{\prime}\Big(\frac{n-\frac 12}{\Delta}\Big)\Big).
\end{equation}
For any 
complex number $\xi$ we have $(\sin(\pi \xi)/(\pi \xi))^2 \ll e^{2\pi |\text{Im }\xi|}/(1+|\xi|^2)$,  and 
further $f(x)\le 4/x^2$ and $|f^{\prime}(x)|\le 8/(|x|(4+x^2))$,  
whence  we deduce that 
\begin{equation}
\label{Gbound} 
|G_\Delta(x+iy) + f(1/\Delta)| \ll  \frac{\Delta^2}{1+|x+iy|} e^{2\pi |y|}.
\end{equation}

We now cull  from Theorem 1.1  of Carneiro and Vaaler \cite{CV} 
various facts about the function $G_\Delta(z)$.  This function is 
entire of exponential type at most $2\pi$, and for real $x$ we have that $G_\Delta(x) \le f_\Delta(x)$.  
We expect that $G_{\Delta}(x)+f(1/\Delta)$ is non-negative for all real $x$, but for our purposes 
a cruder lower bound suffices.  Since $f(x)\ge 0$ and $f^{\prime}(-x)=-f^{\prime}(x)$, 
by pairing the terms $n\ge 1$ with the terms $1-n\le 0$ we obtain from \eqref{Meas3} 
that 
\begin{eqnarray*}
G_{\Delta}(x) + f(1/\Delta) &\ge& \Big(\frac{\cos(\pi x)}{\pi}\Big)^2 \sum_{n=1}^{\infty} \frac{1}{\Delta}
f^{\prime}\Big(\frac{n-\frac12}{\Delta}\Big)\Big( \frac{1}{x-n+\frac 12} - \frac{1}{x+n-\frac 12}\Big) 
\\
&=& \sum_{n=1}^{\infty} \Big(\frac{\sin^2 (\pi (x-n+\frac 12))}{\pi (x^2-(n-\frac 12)^2)} \Big)
\frac{2(n-\frac 12)}{\Delta} f^{\prime}\Big(\frac{n-\frac 12}{\Delta}\Big),
\\
\end{eqnarray*}
and from this we may easily deduce that there is a constant $C$ such that 
\begin{equation} 
\label{cons0} 
-C \frac{\Delta^2}{\Delta^2 +x^2} \le G_{\Delta}(x) + f(1/\Delta) \le f(x/\Delta).   
\end{equation} 

By part (v) of Theorem 1.1 of \cite{CV} we have 
\begin{multline}
\label{cons6}
\int_{-\infty}^{\infty} (f(x/\Delta) - (G_\Delta(x)+f(1/\Delta))) e^{-2\pi i tx} dx 
\\
= \int_{0}^{\infty} \Big( \frac{2\lam}{\lam^2+4\pi^2 t^2} -{\hat L}(\lam,t)\Big) 
\frac{2(1-\cos(2\Delta \lam))}{\lam} d\lam, 
\\
\end{multline}
where ${\hat L}(\lam,t)=0$ if $|t|\ge 1$ and for $|t|\le 1$ we have (see Lemma 3.2 of \cite{CV})
 \begin{equation} 
 \label{cons3}
 {\hat L}(\lam,t) = \frac{(1-|t|)\sinh(\lam/2) \cos(\pi t) + \frac{\lam}{2\pi} |\sin \pi t| \cosh(\lam/2)}{\sinh^2 (\lam/2) + \sin^2\pi t}. 
 \end{equation}
Now $f(x/\Delta)$ is integrable, and we may check that 
\begin{equation*} 
\label{cons5}
\int_{-\infty}^{\infty} f(x/\Delta)e^{-2\pi itx}dx = \int_{0}^{\infty}\frac{2\lam}{\lam^2 +4\pi^2 t^2} 
\frac{2(1-\cos(2\Delta\lam))}{\lam} d\lam, 
\end{equation*} 
so that $G_{\Delta}(x)+f(1/\Delta)$ is also integrable and 
\begin{equation} 
 \label{cons2} 
 \int_{-\infty}^{\infty} (G_\Delta(x)+f(1/\Delta))e^{-2\pi i xt} dx = \int_{0}^{\infty} {\hat L}(\lam,t) \frac{2(1-\cos(2\Delta\lam))}{\lam} d\lam.
 \end{equation} 
   Since 
 $$
 {\hat L}(\lam,t) \ll \frac{1+\lam}{\sinh(\lam/2)} + \frac{\lam}{(\sinh(\lam/2))^2}, 
 $$
 and $(1-\cos(2\Delta \lam))/\lam \ll \min (1/\lam, \Delta^2 \lam)$, we deduce from 
 \eqref{cons3} and \eqref{cons2} that
 \begin{equation} 
 \label{cons4} 
\Big| \int_{-\infty}^{\infty} (G_\Delta(x)+f(1/\Delta)) e^{-2\pi i xt} dx \Big| \ll \Delta. 
\end{equation}

Moreover from \eqref{cons6} and a little calculus we find that 
 \begin{eqnarray} 
 \label{cons1} 
 \int_{-\infty}^{\infty} (f_\Delta(x)-G_\Delta(x)) dx  
 &=& \int_0^{\infty} \Big(\frac 2x -\frac{1}{\sinh(x/2)} \Big) \frac{2(1-\cos(2\Delta x))}{x} dx \nonumber \\
 &= & 2\log 2 -2\log (1+e^{-4\pi \Delta}).
 \end{eqnarray}

We are now in a position to prove Proposition \ref{mainprop}.   We take 
$g_{\Delta}(z) = G_\Delta(z\Delta) + f(1/\Delta)$, so that for real $x$ we have 
$g_{\Delta}(x) = G_\Delta(x\Delta)+ f(1/\Delta) \le f_\Delta(x\Delta)+f(1/\Delta) =f(x)$.  
Since $G_\Delta$ has exponential type at most $2\pi$, we see that $g_{\Delta}$ 
has exponential type at most $2\pi \Delta$.  Further, ${\hat g}_{\Delta}(t)=\Delta^{-1} {\hat G}_{\Delta}(t/\Delta)$.  Thus part (i) of  Proposition \ref{mainprop} follows from \eqref{Gbound} and \eqref{cons0}, 
part (ii) from \eqref{cons3}, \eqref{cons2} and \eqref{cons4}, and part (iii) from \eqref{cons1}.

 \section{Discussion} 
 
 \noindent The estimate \eqref{mainbound} gives a variant of the main Proposition of 
 \cite{Moment} which states that for large $t$, 
 $$
 \log |\zeta(\tfrac 12+it)| \le \text{Re }\sum_{n\le x} \frac{\Lambda(n)}{n^{\frac 12+\frac{\lam}{\log x}+ it}\log n} \frac{\log (x/n)}{\log x} + \frac{(1+\lam)}{2} \frac{\log t}{\log x} + O\Big( \frac{1}{\log x}\Big),
 $$
 where $2\le x\le t^2$, and  $\lam\ge \lam_0=0.4912\ldots$ where $\lam_0$
 denotes the unique positive real number satisfying $e^{-\lam_0} =\lam_0 +\lam_0^2/2$.
 For large $t$ it is difficult to give good estimates for the sum over $n$ above (or in \eqref{mainbound}) 
 and this is the barrier to establishing better estimates for $|\zeta(\tfrac 12+it)|$.  
 However one can study the frequency with which such sums get large, and this 
 information is used in \cite{Moment} to understand the size of moments of $\zeta(\tfrac 12+it)$.

 In light of our work we can view the Proposition in \cite{Moment} as constructing 
 a different minorant of our function $f(x)$.   We start with (for positive $\alpha$ and $x$ real)
 $$
 K_{\Delta}(\alpha,x)=2\pi \int_{-\Delta}^{\Delta} \Big(1-\frac{|t|}{\Delta} \Big)e^{-2\pi \alpha |t|-2\pi i tx} dt 
 = \frac{2\alpha}{\alpha^2+x^2} -\frac{1}{\pi \Delta} \text{Re}\frac{1-e^{-2\pi \Delta(\alpha+ix)}}{(\alpha+ix)^2}.
 $$
 Integrating both sides from $\alpha_0>0$ to $2$ we obtain 
 $$ 
 \int_{\alpha_0}^2 K_{\Delta}(\alpha,x)d\alpha \le \log \frac{4+x^2}{\alpha_0^2 +x^2} 
 +\frac{1}{\pi \Delta} \Big( \frac{2}{4+x^2} -\frac{\alpha_0}{\alpha_0^2+x^2}\Big) 
 + \frac{1}{\pi \Delta} \int_{\alpha_0}^2  \frac{e^{-2\pi \Delta \alpha}}{\alpha_0^2+x^2}d\alpha,
 $$
 and upon rearranging
 $$ 
\log \frac{4+x^2}{\alpha_0^2+x^2}\ge  \int_{\alpha_0}^2 K_{\Delta}(\alpha,x) d\alpha +\frac{1}{\pi \Delta(\alpha_0^2+x^2)}
\Big( {\alpha_0}  -\frac{e^{-2\pi \alpha_0 \Delta}}{2\pi \Delta} 
\Big) - \frac{1}{\pi \Delta} \frac{2}{4+x^2}. 
$$  
 Since $\log ((\alpha_0^2+x^2)/x^2) \ge \alpha_0^2/(\alpha_0^2+x^2)$ we 
 conclude that 
 $$ 
 f(x) \ge \int_{\alpha_0}^{2} K_{\Delta}(\alpha,x) d\alpha +\frac{1}{\alpha_0^2+x^2} 
 \Big(\alpha_0^2 + \frac{\alpha_0}{\pi \Delta} - \frac{e^{-2\pi \alpha_0 \Delta}}{2\pi^2 \Delta^2}\Big) - 
 \frac{1}{\pi \Delta} \frac{2}{4+x^2}.
 $$ 
 If we choose $\alpha_0 \ge \lam_0/(2\pi \Delta)$ then the middle term above is non-negative and 
 we have shown that for such $\alpha_0$ 
 $$ 
 f(x) \ge \int_{\alpha_0}^{2} K_{\Delta}(\alpha,x) d\alpha-\frac{1}{\pi \Delta} \frac{2}{4+x^2}. 
 $$ 
 The first term in the RHS above clearly has Fourier transform supported in $[-\Delta,\Delta]$.  
 The second term may be easily approximated by functions having compactly supported 
 Fourier transform; for example, assuming that $2\pi \Delta \ge 2$ say, we can see 
 from the definition of $K_\Delta$ that $2/(4+x^2) \ge \frac 12 K_{\Delta}(2,x) (1-1/2\pi \Delta)^{-1}$, 
 so that $f(x) \ge \int_{\alpha_0}^2 K_{\Delta}(\alpha, x) d\alpha - K_{\Delta}(2,x)/(2\pi \Delta -1)$.  
 Using the explicit formula with such a minorant gives an alternative proof of the 
 Proposition in \cite{Moment}.  Although the construction of minorants given above (which 
 amounts to taking convolutions with functions whose Fourier transforms have desired 
 compact support) is not optimal, the method works  
 for related functions such as $\log ((4+x^2)/(\alpha^2+x^2))$ (this  arises 
 in bounding $\log |\zeta(\tfrac 12+\alpha+it)|$) which do not fit the framework of Carneiro 
 and Vaaler.

 Theorem \ref{zetabound} may be extended to general $L$-functions.   To be concrete,  
 consider the framework described in Chapter 5 of \cite{IK}.   Thus we consider 
 $L$-functions given in Re $s>1$ by the absolutely convergent series and product
 $$
 L(f,s) = \sum_{n=1}^{\infty} \frac{\lam_f(n)}{n^s} = \prod_{p}\prod_{j=1}^{d} \Big(1-\frac{\alpha_j(p)}{p^s}\Big)^{-1} ,
 $$
   where the `degree' $d$ is a fixed natural number.  We assume that there is an integer 
   $q(f)\ge 1$ and complex numbers $\kappa_j$ with Re$(\kappa_j) >-1$ such that 
   $$
   \Lambda(f,s) = \Big(\frac{q(f)}{\pi^d}\Big)^{\frac s2} \prod_{j=1}^{d} \Gamma\Big(\frac{s+\kappa_j}{2}
   \Big) L(f,s)
   $$
   is entire of order $1$ except possibly for poles at $s=0$ and $1$.  Moreover we suppose 
   that a functional equation 
   $$ 
   \Lambda(f,s) = \epsilon(f) \Lambda({\overline f},1-s), 
   $$ 
   holds, where $\epsilon(f)$ is a complex number of size $1$, and $\Lambda(\overline{f},s) 
   = \overline{\Lambda(f,\overline{s})}$.   We assume the Generalized Riemann Hypothesis 
   for $L(f,s)$, namely that the zeros of $\Lambda(f,s)$ all lie on the line Re$(s)=\frac 12$, 
   and then seek a bound for $L(f,\tfrac 12)$ in 
   terms of the {\sl analytic conductor} $C(f):= q(f) \prod_{j=1}^{d} (3+ |\kappa_j|)$.   
   Making minor modifications to our argument we find that 
   \begin{multline}
   \label{Lbound}
   \log |L(f,\tfrac 12)| \le \frac{\log C(f)}{2\pi \Delta} \log \Big(\frac{2}{1+e^{-4\pi \Delta}}\Big) 
   + \frac{1}{4\pi } \sum_{n=2}^{\infty} \frac{1}{\sqrt{n}} g_{\Delta}\Big(\frac{\log n}{2\pi}\Big) 
   (\Lambda_f(n)+ \Lambda_{\overline{f}}(n))\\
    + O\Big( \frac{\Delta^2 e^{\pi \Delta}}{1+\Delta C(f)} + 1\Big), 
   \end{multline}
   where $\Lambda_f(n)$ and $\Lambda_{\overline{f}}(n)$ are the Dirichlet series coefficients 
   of $-L^{\prime}/L(f,s)$ and $-L^{\prime}/L(\overline{f},s)$ respectively.
   If we now assume the Ramanujan conjectures (which imply that $|\Lambda_f(n)| \le d\Lambda(n)$) 
   then, choosing $\pi \Delta = (1-o(1)) \log \log C(f)$ and estimating the 
   sum over $n$ in \eqref{Lbound} trivially, we obtain that 
   $$
   \log |L(f,\tfrac 12)| \le   \Big( \frac{\log 2}{2} +o(1)\Big) \frac{\log C(f)}{\log \log C(f)}, 
   $$ 
  which is the analog of Theorem \ref{zetabound}.  If we do not assume the Ramanujan conjectures, 
  then using that $|\Lambda_f(n)| \le dn \Lambda(n)$ (which follows from our 
  assumption that the Euler product converges absolutely in $\text{Re}(s)>1$), and choosing $\pi \Delta=(\tfrac 13-o(1)) \log \log C(f)$ we obtain 
  $$ 
  \log |L(f,\tfrac 12)| \le  3 \Big( \frac{\log 2}{2} +o(1)\Big) \frac{\log C(f)}{\log \log C(f)}. 
  $$

\end{document}